\magnification=\magstep1
\input amstex
\documentstyle{amsppt}
\nologo
\NoBlackBoxes
\pageheight{8.5truein}
\pagewidth{6.5truein}

\document

\topmatter
\title
  On the Balazard-Saias criterion for the Riemann Hypothesis 
  \endtitle
\author J. B. Conrey and G. Myerson
\endauthor
\address 
 %$\dag$
J. B. Conrey\endgraf
American Institute of Mathematics\endgraf
360 Portage Ave.\endgraf
Palo Alto, CA 94306\endgraf
{\it E-mail address:}  {\bf conrey\@aimath.org\endgraf}
\null
Department of Mathematics\endgraf
Oklahoma State University\endgraf
Stillwater, OK 74078-0613\endgraf
\null
G. Myerson \endgraf
Department of Mathematics \endgraf
Macquarie University \endgraf
NSW 2109 Australia \endgraf
{\it E-mail address:} {\bf gerry\@mpce.mq.edu.au}
\endaddress

\thanks The visit of the first author to Macquarie University was supported
by a Macquarie University research grant.  The research of the first author is
also supported by the American Institute of Mathematics and the NSF.
\endthanks

\endtopmatter 
\def\starsum{\mathop{{\sum}^*}}

\def\c{\raise.4ex\hbox{$\chi$}}

 \head
 1. Introduction
 \endhead
 
 Recently, Balazard and Saias [BS2] have shown that 
$$
\lim_{N\to \infty} \inf _{D_N}\int _{-\infty}^\infty 
\bigg| \frac{1-\zeta(\tfrac 12+it)D_N(\tfrac 12+it)}
{\tfrac 12+it}\bigg|^2~dt=0 
$$
   implies the Riemann Hypothesis, where 
   $$ 
D_N(s):=\sum_{n\le N}\frac {d_n}{n^s}
$$
   ranges over all Dirichlet polynomials of length $N$.
   
   It is natural that one may wish to investigate this integral 
   taking for $D_N$  
   a partial sum of the Dirichlet series for $1/\zeta(s)$, 
   $$ 
\sum_{n\le N}\frac {\mu(n)}{n^s}.
$$
   However, this choice has some deficiencies, mainly due to the sharp
cutoff of the sum at~$N$, and it is known that this choice does not lead to 
the desired conclusion.
   
   A better choice is $D_N=M_N$ where
   $$
M_N(s):=\sum_{n\le N}\frac {\mu(n)\frac{\log  (N/n)}{\log N}}{n^s}
=\sum_{n\le N}\frac{b_n}{n^s}.
$$
    $M_N$ has its origins in the works of Selberg and is the mollifier used 
in Levinson's work on critical zeros of the Riemann zeta-function.
    Recently, Conrey and Farmer (in preparation) have shown that
    if the Riemann Hypothesis is true and if the zeros of $\zeta(s)$ are 
separated from each other, in the sense that there is a $\delta>0$ such that
    for each zero $\rho$ the derivative of $\zeta$
    satisfies
    $$
|\rho|^{1-\delta}|\zeta'(\rho)| \gg 1,  
$$
    then
    $$
\lim_{N\to \infty}  \int _{-\infty}^\infty 
\left| \frac{1-\zeta(\tfrac12+it)M_N(\tfrac 12+it)}
{\tfrac 12+it}\right|^2~dt =0. 
$$

It is not difficult to deduce by the criterion of Balazard and Saias that the 
Riemann Hypothesis follows from
 $$
   \lim_{N\to \infty} \frac{1}{2\pi} \int _{-\infty}^\infty 
   \left| \frac{\zeta(\tfrac12+it)M_N(\tfrac 12+it)}{\tfrac12+it}\right|^2
~dt =1.
$$
(Square out the integrand and use Cauchy's theorem to evaluate the easy terms 
that arise.)

 \proclaim{Proposition 1} We have
 $$ 
\frac 1{2\pi} \int _{-\infty}^\infty 
   \bigg| \frac{\zeta(\tfrac12+it)M_N(\tfrac12+it)}{\tfrac12+it}\bigg|^2~dt =
\int_0^\infty\bigg| \sum_{n\le N}\frac{b_n \{nu\}}{n} \bigg|^2 
\frac {du}{u^2}
$$
where 
$$
\{x\}=x-[x]
$$ 
is the fractional part of $x$.
\endproclaim
\demo{Proof} The left side is  
$$
\sum_{h, k\le N}\frac{b_h b_k}{k} F(h/k)
$$ 
where
$$
F(x)=\frac{1}{2\pi i}\int_{(\tfrac 12)} \frac{\zeta(s)\zeta(1-s)}{s(1-s)}
x^{-s} ~ds;
$$
the notation $(\tfrac12)$ stands for the vertical path from $\tfrac 12 -i\infty$ to $\tfrac 12+i\infty$. 
Now $F$ may be expressed as a convolution 
$$
F(x)=\int_0^\infty f(u)g(x/u)\frac{du}{u}
$$
where 
$$ 
f(u)=\frac{1}{2\pi i} \int_{(\tfrac1 2)}\frac{\zeta(s)}{s} u^{-s}~ds=
\frac{-1}{u} +\left[\frac 1u\right]
$$
and 
$$ 
g(u)=\frac{1}{2\pi i} \int_{(\tfrac1 2)}\frac{\zeta(1-s)}{1-s} u^{-s}~ds=
\frac{1}{u}(-u+[u]).
$$
By a change of variable
$$
F(h/k)=\frac 1 h \int_0 ^\infty \{hu\}\{ku\}\frac{du}{u^2},
$$
and the proposition follows.  
\enddemo

Thus, it is natural to ask about the series
$$ 
W_N(\alpha)= 
\sum_{n\le N}  \frac {\mu(n)\frac{\log (N/n)}{\log N}\{n\alpha\}}{n}.\tag1
$$
 
In this paper we show in Theorem 1 
that
$$
\lim_{N\to \infty} W_N(\alpha) = \frac{-\sin( 2 \pi \alpha)}{\pi}
$$
uniformly for all real $\alpha$.

We remark that
$$
\int_{0}^\infty \left(\frac{\sin( 2 \pi u)}{\pi }\right)^2 \frac{du}{u^2}=1
$$
but see Remark 1 after Theorem 2.

This research was carried out while the first author 
was visiting Macquarie University. 
He thanks the Department of Mathematics at Macquarie University
for its hospitality during a very pleasant visit. 

\head
2. Heuristics and statements of theorems
\endhead

The series in (1) breaks up   into 
$W_N(\alpha)=U_N(\alpha)-\frac{1}{\log N}V_N(\alpha)$ where 
$$
U_N(\alpha)=\sum_{n\le N} \frac {\mu(n)\{n\alpha\}}{n}
$$
and
$$
V_N(\alpha)= \sum_{n\le N} \frac {\mu(n)\{n\alpha\}\log n}{n}.
$$ 

To motivate our work we observe that
by the prime number theorem,
$$
\sum_{n=1}^\infty \frac{\mu(n)}{n}=0,
$$
so  that
$$ 
U_N(\alpha)=\sum_{n=1}^N \frac {\mu(n)\psi(n\alpha)}{n}+o(1)\tag2
$$
where the saw-tooth function $\psi(x)$ is defined to be zero at integer 
arguments and
$$  
\psi(x)=x-[x]-1/2 =-\sum_{m=1}^\infty \frac{\sin (2\pi m x)}{\pi m} 
$$
 for non-integral $x$.  If we naively insert this series for $\psi(x)$ into 
the sum in (2) and group terms with $mn=k$ we are led to guess that
$$ 
\split
\sum_{n=1}^\infty \frac {\mu(n)\psi(n\alpha)}{n} &=
 -\sum_{n=1}^\infty\frac{\mu(n)}{n}\sum_{m=1}^\infty \frac{\sin (2 \pi m n
\alpha)}{\pi m}\\
 &= -\sum_{k=1}^\infty \frac{\sin (2 \pi k \alpha)}{\pi k}
 \sum_{n\mid k} \mu(k)=-\frac1 \pi \sin( 2\pi \alpha). 
 \endsplit
 $$

 The series involved are only conditionally convergent so that the
interchange of summation is not easily justified. 
 
 In [D1] and [D2], Davenport addressed the question of the convergence of
$U_N(\alpha)$. In the first paper, he showed that 
 $$
\lim_{N\to \infty} U_N(\alpha) =-\frac 1 \pi \sin (2 \pi \alpha) 
$$
 for almost all $\alpha$. In the second paper, after Vinogradov's methods
were developed, he showed that the formula is true for all real $\alpha$ and 
the convergence is uniform. In 1976 S.~Segal [S] showed how to derive the
formula from a Mellin transform. His method does not seem to show that the
convergence is uniform.
  
  A similar argument for 
   $$ 
V^*_N(\alpha):=\sum_{n \le N} \frac{\mu(n)\log n \psi(n\alpha)}{n}
$$
    leads one  to guess that
 $$
\lim_{N\to \infty} V^*_N(\alpha)=
\sum_{n=1}^\infty \frac{\Lambda(n)\sin (2 \pi n \alpha)}{\pi n}.
 $$
 
 Davenport did not address this particular series. Segal's theorem
 is rather general and shows that the identity above holds in the sense 
 that if either side converges, then so does the other side and to
 the same value.
 
 It is the goal of this paper to prove
 
 \proclaim{Theorem 1} Let 
 $$
 W_N(\alpha)=
\sum_{n=1}^N \frac {\mu(n)\frac{\log (N/n)}{\log N}\{n\alpha\}}{n}.
 $$
  Then,
 $$ 
\lim_{N\to \infty} W_N(\alpha) =-\frac {\sin (2 \pi \alpha)}{\pi}
$$
 uniformly for all real $\alpha$.
 \endproclaim

In order to do accomplish this goal, we need the following result, which is of 
independent interest (see Remark 2).
  \proclaim{Theorem 2} The series 
  $$  
T(\alpha)= \sum_{n=1}^\infty \frac{\Lambda(n)\sin (2 \pi n \alpha)}{\pi n}
 $$
  converges for all real $\alpha$.  The convergence is bounded in
  the sense that there is an absolute constant $c>0$ such that the partial
sums
  $$
\bigg|\sum_{n\le N} \frac{\Lambda(n)\sin (2\pi n \alpha)}{\pi n}\bigg|\le c
$$
   for all $N$ and $\alpha$.
   \endproclaim

 {\bf Remark 1.}
 We   cannot conclude that the Riemann Hypothesis holds
   because we cannot show that 
 $$
 \frac 1{u^2}\sum_{n\le N}  \frac {\mu(n)\frac{\log (N/n)}{\log N}\{nu\}}{n}
 \to -\frac{\sin (2\pi u)}{\pi u^2}
 $$
 uniformly. In fact, one can see that if $0< u <1/N$ then
 $$
\sum_{n\le N}  \frac {\mu(n)\frac{\log (N/n)}{\log N}\{nu\}}{n}
 =u\sum_{n\le N} \mu(n) \frac{(\log N/n)}{\log N}
$$
 so that the integral from 0 to $1/N$ of the square of this expression is
just
 $$
\frac{1}{N} \bigg|\sum_{n\le N} \mu(n) \frac{\log (N/n)}{\log N}\bigg|^2.
$$
 The sum over $n$ has an explicit formula; it is
 $$
\frac{1}{2 \pi i \log N} \int_{(c)} \frac{N^s}{\zeta(s)}\frac{ds}{s^2}
 = \frac{1}{\log N}\sum_{\rho} \frac{N^\rho}{\zeta'(\rho) \rho^2} + o(1),
$$
 say, on assuming that the zeros are simple and that 
$|\zeta'(\rho) \rho|\gg |\rho|^ \delta$ for some $\delta>0$ (the integral is
from $c-i\infty$ to $c+i\infty$ where $c>1$).  In this case the series is
absolutely convergent and the size of the sum depends on 
 $\sup_\rho |N^\rho|$. If the Riemann Hypothesis is true, this series is 
bounded uniformly by $N^{1/2}$ from which it follows that
 $$
\sum_{n\le N} \mu(n) \frac{(\log N/n)}{\log N} \ll \frac{N^{1/2}}{\log N}
$$
 and so the integral from 1 to $1/N$ is $\ll 1/\log^2N$. 
 The upshot is that handling the integral over this beginning range clearly
depends on the Riemann Hypothesis.
 
 {\bf Remark 2.} The function $T(\alpha)$ seems to be rather interesting.
It appears to be  continuous at all irrationals, and to have  a jump
discontinuity at $a/q$,   with a jump on either side of size 
$\tfrac12\mu(q)/\phi(q)$  and to satisfy
 $$
T\left(\tfrac a q \right) =\lim_{n\to \infty} \tfrac 1 2 
 \left(T\left(\tfrac a q +\tfrac{1}{n}\right)
+T\left(\tfrac a q -\tfrac{1}{n}\right)\right).
$$
 However, we have not proven these assertions.
 \head
 3. Preliminaries
 \endhead
 
 In Davenport's paper it is remarked that it is easy to use the theory of 
 $L$-functions to show that
 $$
\lim_{N\to \infty} U_N(a/q)=-\frac 1\pi \sin (2 \pi a/q)
$$
 for rational $a/q$.  He does not give the proof. Though it is strictly
speaking not needed for what we do, we believe that it is instructive
nevertheless. Thus, we will show, using the theory of $L$-functions, 
   \proclaim{Proposition 2} If $(a,q)=1$, then   
   $$
\lim_{N\to\infty}U_N(a/q)=\lim_{N\to \infty} 
\sum_{n\le N}\frac{\mu(n)\{na/q\}}{n}=\frac{-\sin (2\pi a/q)}{\pi}
    $$
    \endproclaim
 
 For a Dirichlet character $\c$ modulo $q$ the Dirichlet $L$-function is 
defined for $s=\sigma+it$ with $\sigma >1$ by
  $$
L(s,\c)=\sum_{n=1}^\infty \frac{\c(n)}{n^s}. 
$$
      If $q>1$, then $L(s,\c)$ can be analytically continued as an entire 
function. If $q=1$, then $L(s,\c) =\zeta(s)$ has a simple pole at $s=1$ but is
analytic everywhere else.

   \proclaim{Proposition 3} If $(a,q)=1$,  then 
   $$
   \split 
   & \lim_{N\to \infty}\sum_{n\le N}\frac{\mu(n)\log n \psi(na/q)}{n}
   \\
   &\qquad=
   \frac{1}{\pi i\phi(q)}\sum\Sb \c \bmod q\\ \c\text{odd}
   \endSb
   \c(a)\tau(\overline{\c})\frac{L'}{L}(1,\c)+
   \sum_{p\mid q} \log p\sum_{k=1}^\infty \frac{\sin( 2\pi a p^k/q)}{\pi p^k}
   \endsplit
   $$
   \endproclaim
   
   \proclaim{Proposition 4} If $(a,q)=1$, then
   $$
   \sum_{n=1}^\infty \frac {\Lambda(n) \sin (2\pi a n/q)}{n}=
   \frac{1}{\phi(q)}\sum \Sb \c \bmod q\\ \c\text{ odd} \endSb
   \frac{\c(a)\tau(\overline{\c})}{  i}\frac{L'}{L}(1,\c)+
   \sum_{p\mid q}\log p \sum_{k=1}^\infty \frac{\sin (2\pi a p^k/q)}{p^k}.
   $$
   \endproclaim
   
 {\bf Remark} It is not difficult to give a finite expression for
$\frac{L'}{L}(1,\chi)$, namely 
 $$
\frac{L'}{L}(1, \overline{\chi}) =\log 2\pi +\frac{\gamma}{2} 
 +\frac{\sum_{a=1}^q \chi(a) \log \Gamma (\tfrac a q)}{\sum_{a=1}^q \chi(a) 
\tfrac a q },
$$
where $\gamma\/$ is Euler's constant. 
 
   We also need
   \proclaim{Proposition 5} There is an absolute constant $c_1>0$ such that
   the   sums $V_N(\alpha)$ satisfy
   $$ 
|V_N(\alpha)|\le c_1  
$$
   for all $N\ge 1$ and all $\alpha$.
   \endproclaim

 The basic idea of the proofs of Propositions 2 -- 4 is to use the fact
that $\{na/q\}$ is a periodic function of $n$ with period $q$. We capture the 
 arithmetic progressions modulo divisors of $q$ by using characters, and 
eventually we arrive at an expression involving Dirichlet $L$-functions for odd
characters  at the special values 0 and 1. We make use of the functional
equation for the $L$-function to arrive at the result.

  We can express $L(s,\c)$ in terms of the Hurwitz zeta-function, defined
for $\alpha>0$ and $\sigma>1$ by
  $$
\zeta(s,\alpha) =\sum_{n=0}^\infty \frac{1}{(n+\alpha)^s}.
$$
  The formula is 
  $$
L(s,\c)=q^{-s} \sum_{b=1}^q \c(a) \zeta(s,b/q).
$$
  Since $\zeta(0,b/q)= 1/2-b/q$ (see [WW], section 13.21) we have
  $$
L(0,\c)=\sum_{b=1}^q \c(b)(1/2-b/q).
$$
  
  \proclaim{Lemma 1} Let $\chi$ be a primitive character. Then
  $$
L(0,\overline{\c})L(1,\c)^{-1}=
  \cases \frac{\tau(\overline{\chi})}{\pi i} 
& \text{if ${\raise.4ex\hbox{$\chi$}}$ is odd} \\
  0 & \text{if ${\raise.4ex\hbox{$\chi$}}$ is even}  \endcases
  $$
  \endproclaim
  \demo{Proof}
  If $\c$ is an even primitive character and $q>1$, then
 $$
L(0,\c)=-\frac{1}{q}\sum_{b=1}^q b\c(b) =0.
$$
 If $q=1$, then 
 $$ 
L(1,\c)^{-1}=\zeta(1)^{-1}=0.
$$
 Thus, the formula is true if $\c$ is even.
 
 If $\c$ is an odd primitive character, then $L(s,\c)$ satisfies
  the functional equation (see [D])
  $$ 
\pi^{-\frac12(2-s)}q^{\frac12 (2-s)}
\Gamma\left(\tfrac{2-s}{2}\right) L(1-s,\overline{\c})
  =\frac{iq^{\frac 12}}{\tau(\c)}\pi^{-\frac12 (s+1)}q^{\frac 12 (s+1)}
  \Gamma\left(\tfrac{s+1}{2}\right) L(s,\c) 
$$
  where $\tau(\c)$ is the Gauss sum
  $$
\tau(\c)=\sum_{b=1}^q \c(b) e(b/q)
$$
  with the usual notation $e(x)=e^{2\pi i x}$. We put $s=0$ into this formula,
and use the facts $\Gamma(1/2)=\pi^{1/2}$ and
  $$
\tau(\c) \tau(\overline{\c}) =\c(-1) q
$$
  to obtain the formula in this case.
  \enddemo
  
   \proclaim{Lemma 2} For $(a,q)=1$ we have
   $$
\frac{1}{i\phi(q)} \sum \Sb \c\bmod q\\\c\text{odd}\endSb
   \c(a)\tau(\overline{\c})=\sin( 2\pi a/q).
$$
   \endproclaim
   \demo{Proof} We have
  
   $$
   \split
    \frac{1}{i\phi(q)}\sum\Sb \c\bmod q\\ \c\text{odd} \endSb
   \c(a)\tau(\overline{\c})&=
   \frac 1{2i\phi(q)}\sum \Sb \c\bmod q\endSb
   \left(\c(a)-\c(-a)\right)\tau(\overline{\c})\\
   &=
   \frac 1{2i\phi(q)}\sum \Sb \c\bmod q\endSb
   \left(\c(a)-\c(-a)\right)\sum_{b=1}^q\overline{\c}(b)e(b/q)\\
   &= \frac 1{2i\phi(q)}\sum_{b=1}^q e(b/q)
\sum_{\c\bmod q}\overline{\c}(b)\left(\c(a)-\c(-a)\right)
   \\
   &= \frac{1}{2i}\left(e(a/q)-e(-a/q)\right)=\sin (2 \pi a /q)
  .
  \endsplit
  $$
  \enddemo

   \head
   4. Proofs
   \endhead

   \demo{Proof of Proposition 2}      Let
   $$
U^*_N(\alpha)=\sum_{n\le N}\frac{\mu(n)\psi(n\alpha)}{n}.
$$
   By (2), this is equal to $U_N(a/q)+o(1)$. Then 
   $$
U^*_N(a/q)=\sum_{b=1}^q \psi(ab/q)
\sum \Sb n\le N\\ n\equiv b \bmod q\endSb \frac{\mu(n)}{n}.
$$
   We let $g=(n,q)$. Then
   $$  
U^*_N(a/q)=\sum_{g\mid q}\frac{\mu(g)}{g}
\sum \Sb b=1 \\ (b,q)=1\endSb ^{q/g} 
   \psi\left(\frac{ab}{q/g}\right) 
   \sum \Sb n\le N/g\\ n\equiv b \bmod (q/g) \\(n,g)=1 \endSb
\frac{\mu(n)}{n}.
$$
   Since $(b,q/g)=1$ we can express the congruence condition in the sum
over $n$ by using characters modulo $q/g$. Thus, the sum over $n$ is
   $$
\frac{1}{\phi(q/g)}\sum_{\c \bmod ( q/g)} \overline{\c}(b) \sum
   \Sb n\le N/g \\(n,q/g)=1 \endSb \frac{\mu(n)\c(n)}{n} .
$$
We change variables in the sum over $b$ and replace $b$ by $b\overline{a}$ where
$a\overline{a} \equiv 1 \pmod {q/g}.$ We have 
$$ 
U^*_N(a/q)=\sum_{g\mid q}\frac{\mu(g)}{g}
\frac{1}{\phi(q/g)}\sum_{\c \bmod (q/g)} \c(a) \sum_{b=1}^{q/g}
\overline{\c}(b)
\psi\left(\frac{b}{q/g}\right)  \sum \Sb n\le N/g \\(n,q/g)=1
 \endSb \frac{\mu(n)\c(n)}{n} .  
$$

The sum over $b$ is $-L(0,\overline{\c})$.  Thus, 
$$
U^*_N(a/q)=-\sum_{g\mid q}\frac{\mu(g)}{g}
\frac{1}{\phi(q/g)}\sum _{\c \bmod (q/g)} \c(a) L(0,\overline{\c})
\sum \Sb n\le N/g \\(n,q/g)=1
 \endSb \frac{\mu(n)\c(n)}{n} . \tag3 
$$
 
 Recall that $L(0,\c)=0$ if $\c$ is a non-principal character to an even
modulus. So, we can restrict the sum over $\c$ above to characters that are
either odd or principal.

   The sum over $n$ in (3) is
   $$ 
\split
\sum \Sb n\le N/g \\(n,q/g)=1\endSb \frac{\mu(n)\c(n)}{n}
&=\frac{1}{2\pi i} \int_{(c)}\frac{\prod_{p\mid g}
\left(1-\c(p)/p^s\right)(N/g)^s}{L(s+1,\c)} \frac{ds}{s} \\
&\sim  L(1,\c)^{-1}\prod_{p\mid g} \left(1-\frac{\c(p)}{p}\right)^{-1}
 \endsplit   
$$
by the prime number theorem for arithmetic progressions.
      
 Thus, we now have
 $$
U^*_N(a/q)=-\sum_{g\mid q}\frac{\mu(g)}{g}
\frac{1}{\phi(q/g)}\sum \Sb \c \bmod \frac qg\\ \c \text{odd}\endSb
 \c(a) L(0,\overline{\c})L(1,\c)^{-1} +E_N(a/q)  
$$
 where $E_N(a/q)\to 0$ as $N\to \infty$ for fixed $a$ and $q$. 
 
 To further simplify the main term we use   Lemma 1. But first
 we have to reduce to primitive characters. If $\c \bmod q$ is induced by 
$\c_1\bmod q_1$ where $\c_1$ is primitive, then
 $$
L(s,\c)=L(s,\c_1) \prod_{p\mid (q/q_1)}
\left(1-\frac{\c_1(p)}{p^s}\right).
$$
 Thus, we can write our main term as
 $$ 
\split
  &-\sum_{g\mid q}\frac{\mu(g)}{g}
\frac{1}{\phi(q/g)}\sum_{r\mid (q/g)}\starsum \Sb \c \bmod r\\ \c
\text{odd}\endSb
 \c(a) L(0,\overline{\c})\prod_{p\mid \frac{q}{gr}}(1-\overline{\c}(p))
\\
& \qquad \qquad L(1,\c)^{-1}
\prod_{p\mid \frac{q}{rg}}\left(1-\frac{\c(p)}{p}\right)^{-1}
 \prod \Sb p\mid g \\p\nmid \frac{q}{g} \endSb 
\left(1-\frac{\c(p)}{p}\right)^{-1}
 \endsplit 
$$
 where the * denotes that the sum is for primitive characters.
 We combine two of the products and use Lemma 1 to rewrite the above as 
 $$
-\frac{1}{\pi i} 
 \sum_{g\mid q}\frac{\mu(g)}{g}
\frac{1}{\phi(q/g)}\sum_{r\mid (q/g)}
\starsum \Sb \c \bmod r\\ \c \text{odd}\endSb 
 \c(a) \tau(\overline{\c})
 \prod_{p\mid \frac{q}{gr}}(1-\overline{\c}(p))
 \prod_{p\mid \frac{q}{r}}\left(1-\frac{\c(p)}{p}\right)^{-1}
 .
 $$
 We exchange the orders of summation of $g$ and $r$ and expand one 
 of the products to see that the above is
 $$
 -\frac{1}{\pi i} 
 \sum_{r\mid q}\starsum \Sb \c \bmod r\\ \c \text{odd}\endSb
 \prod_{p\mid \frac qr}\left(1-\frac{\c(p)}{p}\right)^{-1}
 \sum_{d\mid \frac{q}{r}} \mu(d)\overline{\c}(d) \sum_{g\mid \frac{q}{rd}}
 \frac{\mu(g)}{g\phi( q/ g)}
 .
 $$
 
 The sum over $g$ is
 $$ 
\cases \frac{\mu^2\left(\frac{q}{rd}\right)rd}{q\phi(q)} 
&\text{if $(rd,q/rd)=1$}\\
 0 &\text{if $(rd,q/rd)>1$} 
 \endcases
 $$
 
 Thus, the sum over $d$ is
 $$
\sum\Sb d\mid \frac{q}{r}\\ (rd,q/rd)=1\endSb 
 \mu(d)d\overline{\c}(d)\mu^2\left(\frac{q}{rd}\right).
$$
 
 If $(r,q/r)>1$, then this sum is 0 because if $p\mid r$ and $p\mid q/r$,
then $p\mid d$ (since otherwise $p\mid \frac{q}{rd}$), but then $\c(d)=0$ since
 $\c$ is a character modulo $r$. Moreover, the sum is 0 if $q/r$ is not 
squarefree: for if $p^2\mid \frac qr$, then $p^2\mid d $ implies $\mu^2(d)=0$, 
 $p\mid\mid d$ implies $(d,q/rd)>1$, and $p\nmid d$ implies $\mu^2(q/rd)=0$.
 
 Thus, our main term can be rewritten as 
 $$
-\frac{1}{\pi i q \phi(q)} 
 \sum \Sb r\mid q \\ (r,q/r)=1\endSb r\mu^2\left(\frac qr\right)
 \starsum \Sb \c \bmod r\\ \c \text{ odd}\endSb \c(a)\tau(\overline{\c})
 \prod_{p\mid \frac qr}\left(1-\frac{\c(p)}{p}\right)^{-1}
 \prod_{p\mid \frac qr}(1-p\overline{\c}(p)).
$$
 
 Now
 $$
\frac{1-p\overline{\c}(p)}{1-\frac{\c(p)}{p}}
=\frac{p\c(p)-p^2}{p\c(p)-\c(p)^2}=
 -p\overline{\c}(p)
$$
 so that the products over $p$ reduce to
 $$ 
\frac q r\mu\left(\frac q r\right)\overline{\c}\left(\frac q r\right).
$$
 
 Thus,  our main term can now be written as
   $$
-\frac{1}{\pi i   \phi(q)} 
 \sum \Sb r\mid q \\ (r,q/r)=1\endSb  
 \starsum \Sb \c \bmod r\\ \c \text{odd}\endSb \c(a)\mu \left(\frac
qr\right)
 \overline{\c}\left(\frac qr\right)\tau(\overline{\c})
 .
 $$
 Now if $\c \bmod q$ is induced by $\c_1 \bmod r$ then $\tau(\c)=0$ if
 $(r,q/r)>1$ or if $\mu(q/r)~=~0$. If $(r,q/r)=1$ and $q/r$ is squarefree,
 then 
 $$
\tau(\c)=\mu \left(\frac qr\right)
  \c_1\left(\frac qr\right)\tau( \c_1).
$$
  
  Thus, the above expression for our main term simplifies to
  $$
-\frac{1}{\pi i \phi(q)}\sum \Sb \c \bmod q\\ \c \text{odd} 
  \endSb \c(a) \tau(\overline{\c}).
$$
  Now
  $$
\sum \Sb \c \bmod q\\ \c \text{odd} 
  \endSb \c(a) \tau(\overline{\c})
  =\frac 1 2\sum_{\c \bmod q} (\c(a)-\c(-a))\tau(\overline{\c}).
$$
  Also,
  $$ 
\sum_{\c \bmod q}  \c(a) \tau(\overline{\c})
  = \sum_{\c \bmod q} \sum_{b=1}^q \overline{\c}(b) e(b/q)
  =\phi(q) e(a/q).
$$
  
  Thus, 
  the main term reduces to
  $$
\frac{-\sin (2\pi a/q)}{\pi} 
$$
  as desired.
  
\enddemo

 \demo{Proof of Proposition 3} We  reduce this Proposition to several
  instances of Proposition 2.  To do this, we write
     $$
\split
   W_N(a/q) 
&=\sum_{n\le N}\frac{\mu(n)\log\big(\frac{n}{(n,q)}\big)\psi\left(
 {na}{q}\right)}{n}+
  \sum_{n\le N}\frac{\mu(n)\log (n,q)\psi\left(
 {na}{q}\right)}{n}\\
 &= \Sigma_1+\Sigma_2
 \endsplit
 $$
 say. We handle $\Sigma_1$ much as in the proof of Proposition 2.
  We split the range of summation into arithmetic progressions $b \bmod q$
and split further according to the greatest common divisor $g=(b,q)=(n,q)$.
 Thus, we arrive at
 $$
\Sigma_1= 
  -\sum_{g\mid q}\frac{\mu(g)}{g}
\frac{1}{\phi(q/g)}\sum _{\c \bmod \frac qg} \c(a) L(0,\overline{\c})
\sum \Sb n\le N/g \\(n,q/g)=1
 \endSb \frac{\mu(n)\c(n)\log n}{n} .  
$$
 
 Now
 $$
 \split \sum\Sb n=1\\(n,g)=1\endSb ^\infty \frac{\mu(n)\c(n)\log n}{n} 
 &= \left.\frac{d}{ds}L(s,\c\c_{0,g})^{-1}\right|_{s=1}
 \\
 &=  
  L(1,\c)^{-1}\prod_{p\mid g}\left(1-\frac{\c(p)}{p}\right)^{-1}
 \frac{L'}{L}(1,\c\c_{0,g}), 
 \endsplit
   $$
where $\c_{0,g}$ is the principal character modulo~$g$. We can replace the sum
over $n$ with this expression and have exactly the same error term $E_N(q)$ as
in Proposition 2.
 
 We reduce to primitive characters and use Lemma 1, much as before. The
main term of $\Sigma_1$ is then
 $$
  -  \frac{1}{\pi i} 
 \sum_{r\mid q}\starsum \Sb \c \bmod r\\ \c \text{odd}\endSb
 \prod_{p\mid \frac qr}\left(1-\frac{\c(p)}{p}\right)^{-1}
 \sum_{d\mid \frac{q}{r}} \mu(d)\overline{\c}(d) \sum_{g\mid \frac{q}{rd}}
 \frac{\mu(g)}{g\phi\left(\frac q g\right)}
  \frac{L'}{L}(1,\c\c_{0,g})
 $$
  This term 
can now be treated exactly as in the proof of Proposition 2.
It leads to a contribution of
$$
\frac{1}{\pi i}\sum_{\c~\bmod q}\tau(\overline{\c}) \c(a)
\frac{L'}{L}(1,\c).
$$

 To treat $\Sigma_2$ we use the formula 
 $$
\log n =\sum_{s\mid n} \Lambda(s).
$$
 Thus,
 $$
 \split
 \Sigma_2&= \sum_{n\le N}\frac{\mu(n)\log (n,q)\psi\left(
 {na}{q}\right)}{n}\\
 &=  \sum_{n\le N}\frac{\mu(n) \psi\left(
 {na}{q}\right)}{n}\sum\Sb s\mid q\\s\mid n\endSb \Lambda(s)
 \\
 &= \sum_{s\mid q} \Lambda(s) \sum_{n\le N/s}\frac{\mu(sn) \psi\left(
 {sna}{q}\right)}{sn}.
 \endsplit
 $$
 Clearly, $s$ must be a prime divisor of $q$. We change $s$ to $p$
 and have
 $$
\Sigma_2= -\sum_{p\mid q} \log p \sum\Sb n\le N/p\\
 p\nmid n \endSb \frac{\mu(n) \psi\left(
 {nap}{q}\right)}{pn}.
 $$
  
 Now, for any positive integer $k$ let 
 $$
 r(k)=
 \sum \Sb n\le x \\p\nmid n\endSb
 \frac{\mu(n)}{np^k}\psi\left(\frac{anp^k}{q}\right). 
$$
 Then,
 $$
 \split
 r(k)&=
 \sum \Sb n\le x  \endSb
 \frac{\mu(n)}{np^k}\psi\left(\frac{anp^k}{q}\right)
 -\sum \Sb n\le x \\p\mid n\endSb
 \frac{\mu(n)}{np^k}\psi\left(\frac{anp^k}{q}\right)\\
 &=-\frac 1 \pi \frac{\sin (2\pi a p^k/q)}{ p^k} +o\left(\frac{1}{p^k}\right)+
 \sum \Sb n\le x \\p\nmid n\endSb
 \frac{\mu(n)}{np^{k+1}}\psi\left(\frac{anp^{k+1}}{q}\right)
 \\
 &=-\frac 1 \pi\frac{\sin (2\pi a p^k/q)}{ p^k} +o(1)+r(k+1).
 \endsplit
 $$
 If we apply this relation repeatedly, we end up with
 $$
\Sigma_2= \sum_{p\mid q}\log p\sum_{k=1}^\infty 
\frac{\sin (2\pi a p^k/q)}{\pi p^k}+o(1).
$$
 
 Thus, we have proved  Proposition 3.
 \enddemo
 
 \demo{Proof of Proposition 4}
 We have
 $$ 
 \sum_{n=1}^N \frac{\Lambda(n) \sin( 2\pi a n/q)}{n}=
 \sum \Sb n\le N \\ (n,q)=1\endSb \frac{\Lambda(n) \sin( 2\pi a n/q)}{n}+
 \sum\Sb n\le N \\(n,q)>1\endSb \frac{\Lambda(n) \sin( 2\pi a n/q)}{n}  \tag 4
 $$  
 and 
 $$
 \split \sum \Sb n\le N \\ (n,q)=1\endSb \frac{\Lambda(n) \sin( 2\pi a n/q)}{n} &= \frac{1}{i\phi(q)}\sum \Sb \c \bmod q\\ \c \text{ odd} \endSb
 \c(a) \tau(\overline{\c}) \sum_{n \le N} \Lambda(n)\c(n) \\
 &=\frac{1}{i\phi(q)}\sum \Sb \c \bmod q\\ \c \text{ odd} \endSb
 \c(a) \tau(\overline{\c})\frac{L'}{L}(1,\c)+o(1)
 .
 \endsplit
 $$
 To evaluate the second sum on the right side of (4) we observe that since 
$\Lambda$ is supported on prime powers, it must be the case that $(n,q)$ is a
power of a prime $p$, or else the sum is~0. Thus, we can group the terms
according to primes $p$ dividing $q$. For a given $p$ dividing~$q$ 
 the $n$ for which $p\mid (n,q)$ and $\Lambda(n)\ne 0$ are just $n=p^k$ for
 some $k\ge 1.$ Therefore, the second sum is 
 $$
 \sum_{p\mid q} \log p\sum_{p^k\le N} \frac{\sin (2 \pi a
p^k)}{p^k}\sim 
 \sum_{p\mid q} \log p\sum_{k=1}^\infty \frac{\sin( 2 \pi a p^k)}{p^k} 
  .
 $$
  \enddemo
  
  To prove Proposition 5 we use the ideas of   Davenport [D1] and [D2] .
    First, we prove 
  \proclaim{Lemma 3}
  We have
  $$
\bigg|\sum \Sb n \le N\\ q\mid n\endSb \frac{\mu(n)\log n}{n}\bigg|\ll
   \cases \frac 1 {\phi(q)} &\text{if $q\ll \log^hN$ }\\
    \frac{\log N}{q}  &\text{if $q\ge \log^h N$}
    \endcases
$$
    \endproclaim
    
    \demo{Proof}

  To prove this, note that the sum is
  $$
\split
  &\frac{1}{q}\bigg|
\sum \Sb n \le \frac{N}{q}\\ (q,n)=1\endSb \frac{\mu(n)\log nq}{n}\bigg|
  \\
  &\le \frac{1}{q}\bigg|
\sum \Sb n \le \frac{N}{q}\\ (q,n)=1\endSb \frac{\mu(n)\log q}{n}
  \bigg|
  +\frac{1}{q}\bigg|
\sum \Sb n \le N\\ (q,n)=1\endSb \frac{\mu(n)\log n}{n}\bigg|
  \endsplit
  $$
  
  The first term is $O\bigl( (\log q)/q\bigr)$ for all $q$ by [D1] Lemma 1
  and is  $ O\bigl( (\log N)^{-h}\bigr)$ by Lemma 12 of [D2] for $q\le \log^hN$. 
   So it suffices to bound
  $$
\sum \Sb n \le x\\ (q,n)=1\endSb \frac{\mu(n)\log n}{n} .
$$
  Note that 
  $$
\sum \Sb d\mid n\\(d,q)=1\endSb \mu(d)\log d 
=-\Lambda\left(\frac{n}{n_q}\right)
$$
  where $n_q$ is that part of $n$ which is coprime to $q$, i.e., 
$n_q=\prod_{p^k\parallel n,p\nmid q}p^k$. The $n$ for which
$\Lambda(n/n_q)\neq0$   are those of the form $n=dm$ where $d\in q^\infty$ and
$\Lambda(m)\ne 0$ (where $q^\infty$ is the set of all integers all of whose
prime factors divide~$q$).  Thus, 
  $$
  \split
  \sum_{n \le x} \sum\Sb d\mid n\\ (d,q)=1\endSb\mu(d)\log d
  &=-\sum_{d\in q^\infty} \sum_{n\le \frac{x}{d}}\Lambda(n)
  \\
  &\ll
\sum_{d\in q^\infty}\frac{x}{d}\ll x\prod_{p\mid q} \left(1+\frac
1p+\frac1{p^2}+\dots
  \right) \\
  &= x \prod_{p\mid q}(1+1/p) \ll x\prod_{p\le q} (1+1/p)\ll x \log q.
  \endsplit
  $$
  Therefore,
  $$
\bigg| \sum_{n\le x} \sum\Sb d\mid n\\(d,q)=1\endSb
\mu(d)\log d\bigg|\ll x\log q.
$$
  But the left side of this inequality is 
  $$
  \split
  \bigg|
\sum\Sb d\le x\\ (d,q)=1\endSb \mu(d)\log d\left[\frac{x}{d}\right]\bigg|
  &=
  x\bigg|
\sum\Sb d\le x\\ (d,q)=1\endSb \frac{\mu(d)\log d}{d}\bigg|+O(x\log x)
  \\
  &\ll x\log qx.
  \endsplit
  $$
  For $q\le \log^hN$,
  $$
\sum \Sb n\le x \\(n,q)=1\endSb \frac{\mu(n)\log n}{n} 
  = \frac{1}{2\pi i}\int_{(c)} \sum\Sb n =1\\(n,q)=1\endSb^\infty 
  \frac{\mu(n)\log n}{n^s} \frac{x^s}{s}~ds.
  $$
  The series under the integral sign is
  $$
-\frac{d}{ds}\left(\zeta(s)^{-1}\prod_{p\mid q}
\left(1-\frac{1}{p^s}\right)^{-1}\right),
$$
  which, by standard arguments, is
  $$
-\prod_{p\mid q} \left(1-\frac{1}{p}\right)^{-1}
+O\left((\log x)^{-h}\right)
  =-\frac{q}{\phi(q)}+O\left((\log x)^{-h}\right).
$$
  \enddemo

  \proclaim{Lemma 4}
  Let
  $$
V^*_N(\alpha)=\sum_{n\le N} \frac{\mu(n)\log n\psi(n\alpha)}{n}.
$$
  Then for all $N$, $\alpha_1$, $\alpha_2$,
  $$
  |V^*_N(\alpha_1)-V^*_N(\alpha_2)|\ll N\log N |\alpha_1-\alpha_2| +1.
$$
   \endproclaim
   
   \demo{Proof}
   The proof follows Lemma 2 of [D1] as well as Lemmas 12 and 13 of [D2].
   We have that $V_N^*(\alpha)$ is continuous and differentiable, with 
derivative
   $$
\sum_{n \le N} \mu(n)\log n \ll N\log N
$$
   except at rationals $a/q$ with $q\le N$ where it has a jump
discontinuity of size
   $$
-\sum\Sb n\le N\\ q\mid n\endSb \frac{\mu(n)\log n}{n}.
$$
   Thus,
   $$ 
V_N^*(\alpha)-V_N^*(\beta)\ll (\alpha-\beta)N\log N+
   \bigg|\sum_{\alpha\le \frac{a}{q}\le \beta} \sum\Sb n \le N\\q\mid
n\endSb \frac{\mu(n)\log n}{n}\bigg|.
$$
   Now we use the estimates of Lemma 3 for the inner sum and the arguments
   of Lemma 2 of [D1] and Lemma 13 of [D2] to complete the proof.
   \enddemo
   
   \demo{Proof of Proposition 5} Here we follow the proofs of Lemma 14 and 
Theorem 2 of [D2]. Let
   $$ 
   \split
   R_N(\alpha)&=V_N^*(\alpha)-T(\alpha)\\
   &=\sum_{n>N} \frac{\mu(n)\log n\psi(n\alpha)}{n} .
   \endsplit
   $$
   Then  
   $$\int_{\alpha_1}^{\alpha_2} R_N(\alpha) ~ d\alpha =
   \sum_{n>N} \frac{\mu(n)\log n\psi_2(n\alpha_2)}{n}-
   \sum_{n>N} \frac{\mu(n)\log n\psi_2(n\alpha_1)}{n} 
   $$
   where 
   $$ \psi_2(t)=\frac{1}{2\pi^2}\sum_{m=1}^\infty 
   \frac {\cos 2 \pi m t}{m^2}=\int_0^t\psi(u) ~du +\frac{1}{12}.$$
   Thus,  
   $$
   \split
   \sum_{n>N} \frac{\mu(n)\log n\psi_2(n\alpha)}{n}&=
   \frac{1}{2\pi ^2}\sum_{n >N} \frac{\mu(n) \log n}{n^2} 
   \sum_{m=1}^\infty \frac{\cos 2\pi m n \alpha}{m^2} \\
   &=
   \frac{1}{2\pi ^2}\sum_{m=1}^\infty \frac{1}{m^2}
   \sum_{n>N}\frac{\mu(n) \log n \cos 2\pi m n \alpha}{n^2}  \\
   & \ll N^{-1}(\log N)^{-h}
   \endsplit
   $$
   by Theorem 1 of [D2] and partial summation.  Next,
   $$
   \split
   (\alpha_1-\alpha_2)R_N(\alpha_1)
&= \int_{\alpha_1}^{\alpha_2}R_N(\alpha_1) ~d\alpha\\
   &= \int_{\alpha_1}^{\alpha_2}R_N(\alpha) 
   ~d\alpha+\int_{\alpha_1}^{\alpha_2}(R_N(\alpha_1)-R_N(\alpha)) ~d\alpha
   .
   \endsplit
   $$
   Therefore,
   $$
|R_N(\alpha_1)|\le \frac{1}{\alpha_1-\alpha_2}\frac{1}{N\log ^h N}
   +\max_{\alpha_1\le \alpha, \beta\le \alpha_2}
\left|R_N(\alpha)-R_N(\beta)\right|.
   $$
   Now
   $$ 
R_N(\alpha)-R_N(\beta) =V_N^*(\alpha)-V_N^*(\beta)+T(\alpha)-T(\beta)
   \ll 1+ \left| V_N^*(\alpha)-V_N^*(\beta)\right|
$$
   by Theorem 2 of [D2]. Take
   $$ 
\alpha_1-\alpha_2=\frac{1}{N\log^h N}
$$
   and use Lemma 4 to obtain the result.
   \enddemo

  \demo{Proof of Theorem 2}
  Let
  $$
S_u(\alpha) =\sum_{n\le u}\Lambda(n) \sin( 2\pi n \alpha).
$$
  Then
  $$
\sum_{n\le N}\frac{\Lambda(n) \sin( 2 \pi n
\alpha)}{n}=\frac{S_N(\alpha)}{N}
  +\int_2^N\frac{S_u(\alpha)}{u^2}~du.  \tag 5
$$
  Note that 
  $$
|S_N(\alpha)| \le \sum_{n\le N}\Lambda(n)\ll N
$$
  so that the first term on the right side of (5) is uniformly bounded.
  Now let $H>10$ be fixed. Define 
  $$ 
\tau=\tau(u) =\frac u {\log^H u}
$$
  for $u\ge 2$.
  Let $q\le \tau$ be such that
  $$
\left|\alpha-\frac{a}{q}\right| \le\frac{1}{q\tau}
$$
for some $a$.
  Note that for each $u$ there is a unique such $q$. 
  We split the $u $ with $2\le u\le N$ into two sets $R_1(N)$ and $R_2(N)$ 
according to the size of $q$. 
  If $q \le  \log^H u$ then $u\in R_1(N)$, and if
    $ \log^H u \le q \le \tau(u)$, then $u \in R_2(N)$.  We will show that 
    $$
\int_{R_j(N)} \frac{S_u(\alpha)}{u^2}~du 
$$
    is uniformly bounded and has a limit as $N\to \infty$ for $j=1$ and 2.
    
    Suppose $u\in R_2$.  Then, by the theorem of section 25 of [D],
  $$
  \split 
  S_u(\alpha)&\ll \left(\frac u q+u^{4/5}+(uq)^{1/2}\right)\log^4 u \\
  & \ll \frac{u}{(\log u)^{\frac{H}{2}-4}}.
  \endsplit
  $$
  Therefore, 
  $$
\int_{R_2(N)} \frac{S_u(\alpha)}{u^2}~du
  \ll \int_2^N \frac{du}{u\log^{1+\delta}u} \ll 1
  $$
  uniformly for all $N$. The integral over $R_2=\lim_{N\to \infty}R_2(N)$
  is absolutely convergent.
  
  Now suppose that $u\in R_1(N)$. Write
  $$
\alpha=\frac{a}{q}+\beta.
$$
  Then by section 26 of [D],
  $$ 
S_u(\alpha)=\Im \frac{\mu(q)}{\phi(q)}\sum_{n \le u} e(n \beta)
   +O\left(u\exp(-C\sqrt{\log u})\right)
$$
   for an absolute constant $C>0$, where $\Im z\/$ is the imaginary part 
of~$z$. Clearly, the integral over $R_2$ of the  $O$-term   is uniformly
   bounded and converges absolutely.

   Now
   $$
\Im \sum_{n \le u} e(n\beta)=\sum_{n=1}^{[u]} \sin( 2 \pi n \beta)
   =\frac{\sin\big(\tfrac{([u]+1)\beta}{2}\big)\sin 
\big(\tfrac{[u]\beta}{2})}{\sin \beta}.
   $$   
   Thus, for any particular $q$ the integral over $R_2$ of the 
   contribution from the main term above
   is
   bounded by
   $$
   \frac{1}{\phi(q)} \int_2^\infty 
   \left|\frac{\sin\left(\tfrac{([u]+1)\beta}{2}\right)
\sin \left(\tfrac{[u]\beta}{2}\right)}{\sin \beta}\right|
    \frac{du}{u^2}
   .  \tag 6
   $$
  Observe that
  $$
|\sin ([u]\beta) -\sin (u \beta)| \le |\beta|
$$
  and $$
\frac{|\beta|}{|\sin \beta|}\ll1
$$
  so that the expression in (6) is
  $$
  \frac{1}{\phi(q)}\left( \int_2^\infty 
   \left|\frac{\sin^2\left(\tfrac{u\beta}{2}\right)}{\sin \beta}\right|
    \frac{du}{u^2}+O(1)\right).
$$
   Let  $v=u\beta$ to see that the above is bounded by
    $$
  \frac{1}{\phi(q)}\left( \int_2^\infty 
   \left|\frac{\sin^2\left(\frac{u }{2}\right)}{u^2} \right|
    ~du+O(1)\right)\ll\frac{1}{\phi(q)}.
    $$
  All of the $q$ which appear in the above proof are denominators of
  convergents of the continued fraction of $\alpha$. It is easy to see that if
the convergents of $\alpha$ are $p_m/q_m$ then
  $$
\sup_{\alpha} \sum_{n=1}^\infty \frac{1}{\phi(q_m)}\ll 1.
$$
  Thus, the contribution of this part is uniformly bounded and converges.
  
  Thus, we have completed the proof that the partial sums $V^*(N)$ are 
  uniformly bounded.
  
  It only remains to observe that $\lim_{N\to \infty}S_N(\alpha)/N=0$
 for all fixed $\alpha$ to complete the proof of convergence. If $\alpha$ is
rational then convergence of  $T(\alpha)$ follows from Proposition 4. If
$\alpha$  is irrational, then we argue again according to whether 
$N~\in~R_1(N)$ or $N\in R_2(N)$. In the first case, the relevant $q\to \infty$,
and   the second case is clear. Thus, we have convergence in all cases.
  \enddemo
  
  \demo{Proof of Theorem 1}
  It   follows from [D2] that 
  $$
U_N(\alpha)\to \frac{\sin( 2 \pi \alpha)}{-\pi}
$$
  uniformly. Thus, it suffices to show that
  $$
\frac{1}{\log N} V_N(\alpha) \to 0
$$
  uniformly. Hence, it suffices to show that $V_N(\alpha)$
  is uniformly bounded. Proposition~5 shows that $V_N^*(\alpha)$
  is uniformly bounded. If $\alpha$ is irrational then
$$ 
\split
V_N(\alpha)& = V_N^*(\alpha) -\frac{1}{2} 
\sum_{n\le N} \frac{\mu(n)\log n}{n}\\
&=V_N^*(\alpha)-1 +o(1).
\endsplit
$$
If $\alpha=a/q$ is rational, then
$$
 V_N(\alpha) = V_N^*(\alpha) -\frac{1}{2} \sum_{n\le N} 
\frac{\mu(n)\log n}{n}
+\frac{1}{2}\sum \Sb n\le N \\q\mid n\endSb \frac{\mu(n)\log n}{n}.
$$
The last term is uniformly bounded by Lemma 3. 
Thus, $V_N(\alpha)$ is uniformly bounded and the Theorem follows.

  \enddemo

 \enddocument

  \end
  \end{document}